\documentclass[11pt]{amsart}

\usepackage{geometry}                
\usepackage{graphicx}
\usepackage{amssymb}
\usepackage{epstopdf}
\usepackage{tikz}
\usetikzlibrary{calc,shapes.geometric,decorations.markings,shadows,decorations.fractals}
\usepackage{amsaddr}

\newtheorem{thm}{Theorem}[section]

\newtheorem{lemma}[thm]{Lemma}

\theoremstyle{definition}

\title{Factor groups of knot and LOT groups}
\author{R. Gerecke, J. Harlander$^1$, R. Manheimer, B. Oakley, S. Rahman}

\begin{document}
\maketitle
\begin{abstract}A classical result of H. S. M. Coxeter asserts that a certain quotient $B(m,n)$ of the braid group $B(m)$ on $m$ strands is finite if and only if $(m,n)$ corresponds to the type of one of the five Platonic solids. If ${\bf k}$ is a knot or virtual knot, one can study similar quotients $G({\bf k}, n)$ for the corresponding knot group. We identify a class of long virtual knots ${\bf k}$ for which $G({\bf k}, n)$ is infinite for $n\ge 2$. The main feature of these long virtual knots is that their Wirtinger complexes are non-positively curved squared complexes.
\end{abstract}

\section{Introduction}

\let\thefootnote\relax\footnote{$^1$Corresponding author. Contact information can be found at the end of the article.}

This article takes its inspiration from the following beautiful result of H. S. M. Coxeter \cite{Coxeter}, also contained in the collection \cite{Coxeter2}. See also the book ``A Study of Braids" by K. Murasugi and B. I. Kurpita \cite{Braids}, page 81.

\begin{thm} Let $B(m)$ be the braid group on $m$ strands and let $B(m,n)$ be the quotient obtained from $B(m)$ by setting $\sigma^n=1$, where $\sigma$ is any of the standard generators and $n\ge 3$. Then $B(m,n)$ is finite if and only if $(m,n)$ corresponds to the type of one of the five Platonic solids.
\end{thm}

Note that $B(m,1)$ is the trivial group and $B(m,2)$ is the symmetric group $S_m$. In case $m=3$ we have $B(3,3)\cong SL(2,3)$ (type of the tetrahedron), $B(3,4)\cong SL(2,3)\rtimes \mathbb Z_4$ (type of the octahedron), $B(3,5)\cong SL(2,5)\times \mathbb Z_5$ (type of the icosahedron), and $B(3,n)$ is infinite for $n\ge 6$.\\

In this paper we study similar quotients of knot groups. Let $\bf k$ be a knot and let $G({\bf k})$ be the knot group. Let $G({\bf k}, n)$ be the quotient obtained from $G(\bf k)$ by setting $x^n=1$, where $x$ is a meridian of the knot. The series $\{ G({\bf k}, n)\}_{n\in \mathbb N}$ is a knot invariant.
\begin{itemize}
\item Does there exist an $N$ so that $G({\bf k}, n)$ is finite for $n\le N$ and infinite for all other $n$?
Note that in case $\bf k$ is the trefoil knot we have $G({\bf k})\cong B(3)$ and $G({\bf k}, n) \cong B(3,n)$. Hence $G({\bf k},n)$ is finite for $n\le 5$ and infinite for $n\ge 6$. 
\item Which finite groups can occur as $G({\bf k}, n)$? All such finite groups have balanced presentations and hence are {\em interesting} in the sense of Johnson \cite{JohnsonBook}, Chapter 7. As far as we know, no interesting finite groups that require more than three generators are known.
\end{itemize}
We believe these to be difficult questions. In this paper we identify a class of  long virtual knots $\bf k$ for which all $G({\bf k}, n)$, $n\ge 2$, are infinite. An important feature of these virtual knots is the fact that their Wirtinger complexes are non-positively curved square complexes, a feature that classical knots do not have. A convenient way to record Wirtinger complexes is via labeled oriented graphs.\\

A {\em labeled oriented graph} (LOG) is an oriented graph $\mathcal G$ on vertices $\{ a, b, c, ...\}$, where
each oriented edge is labeled by a vertex. All labeled oriented graphs considered in this article are assumed to be {\em compressed}: an edge is not labeled by one of its vertices. Associated with it is the {\em LOG-presentation} 
$$P({\mathcal G})=\langle a, b, c,... \ |\  r_e=1,\ \mbox{$e\in$ edges of $\mathcal G$} \rangle $$ where $r_e=ab(bc)^{-1}$ in case $e=(a|b|c)$ is a labeled edge with initial vertex $a$, terminal vertex $c$, and label $b$. A {\em LOG-complex} 
$K({\mathcal G})$ is the standard 2-complex associated with the LOG-presentation $P(\mathcal G)$, and a {\em LOG-group} $G({\mathcal G})$, is the group defined by the LOG-presentation. A labeled oriented graph is called {\em injective} if each vertex label occurs at most once as an edge label.
A {\em labeled oriented tree} (LOT) is a labeled oriented graph where the underlying graph is a tree. A {\em labeled oriented interval} (LOI) is a labeled oriented graph where the underlying graph is a subdivided interval.\\

There is a one-to-one correspondence between Wirtinger complexes of long virtual knots and LOI-complexes. See Harlander, Rosebrock \cite{HarlanderRosebrock} and the references therein. We will state our results in terms of LOT-complexes.\\

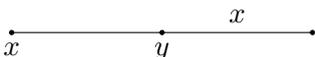
\begin{figure}[here] 
   \centering
\begin{tikzpicture}
\fill (0,0) circle (1pt);
\fill (2,0) circle (1pt);
\fill (4,0) circle (1pt);
\draw (0,0) to (2,0);
\draw (2,0) to (4,0);
\node [below] at (0,0) {$x$};
\node [below] at (2,0) {$y$};
\node [above] at (3,0) {$x$};
\end{tikzpicture}
 \caption{An additional edge combination that we need to avoid. We want an edge labeled $x$ not to be connected to an edge with vertex $x$. }
   \label{fig:npc+}
\end{figure}

A 2-complex $K$ is a {\em non-positively curved squared complex} if the boundary of every 2-cells has four edges and cycles in vertex links have length at least four. This is equivalent to saying that in a reduced surface diagram over $K$ we have at least four 2-cells grouped around every vertex. We will discuss surface diagrams over 2-complexes in more detail in the next section. It was shown by Rosebrock \cite{Rosebrock} that a LOG-complex, associated with an injective labeled oriented graph $\mathcal G$, is a non-positively curved squared complex if and only if the edge combinations shown in Figure \ref{fig:2-cycles} and in Figure \ref{fig:3-cycles} do not occur in $\mathcal G$. Note that the edges in these figures can carry any orientation. Thus Figure \ref{fig:2-cycles} gives rise to four oriented edge combinations, and Figure \ref{fig:3-cycles} gives rise to eight oriented edge combinations.
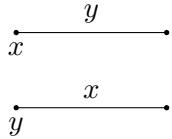
\begin{figure}[here] 
   \centering
   \begin{tikzpicture}
\fill (0,1) circle (1pt);
\fill (2,1) circle (1pt);
\fill (0,0) circle (1pt);
\fill (2,0) circle (1pt);

\draw (0,1) to (2,1);
\draw (0,0) to (2,0);

\node [below] at (0,1) {$x$};
\node [above] at (1,1) {$y$};
\node [below] at (0,0) {$y$};
\node [above] at (1,0) {$x$};
\end{tikzpicture}
\caption{An edge combination that implies 2-cycles in the vertex link of $K({\mathcal G})$. Any edge orientation is possible.}
\label{fig:2-cycles}
\end{figure}
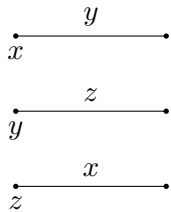
\begin{figure}[htbp] 
   \centering
   \begin{tikzpicture}
\fill (0,0) circle (1pt);
\fill (2,0) circle (1pt);
\fill (0,1) circle (1pt);
\fill (2,1) circle (1pt);
\fill (0,2) circle (1pt);
\fill (2,2) circle (1pt);

\draw (0,0) to (2,0);
\draw (0,1) to (2,1);
\draw (0,2) to (2,2);
\node [below] at (0,0) {$z$};
\node [below] at (0,1) {$y$};
\node [below] at (0,2) {$x$};
\node [above] at (1,0) {$x$};
\node [above] at (1,1) {$z$};
\node [above] at (1,2) {$y$};
\end{tikzpicture}
 
 \caption{An edge combination that implies 3-cycles in the vertex link of $K({\mathcal G})$. Any edge orientation is possible.}
   \label{fig:3-cycles}
\end{figure}
One more edge combination that we need to pay attention to is shown in Figure \ref{fig:npc+}.

Given a labeled oriented graph $\mathcal G$ and $x$ a vertex we denote by $$P({\mathcal G}, x^n)=\langle a, b, c, ...\ |\ x^n=1, r_e=1,\ \mbox{$e\in$ edges of $\mathcal G$} \rangle, $$ by $K({\mathcal G}, x^n)$ the standard 2-complex constructed from $P({\mathcal G}, x^n)$,  and by $G({\mathcal G}, x^n)$ the group defined by $P({\mathcal G}, x^n)$. Note that in case $\mathcal G$ is connected, all vertices give conjugate elements in $G({\mathcal G})$ and hence all $G({\mathcal G}, x^n)$, $x$ a vertex in $\mathcal G$, are identical. We can now state our main result.

\begin{thm}\label{thm:mainthm} Let $\mathcal T$ be an injective labeled oriented tree that does not contain the edge combinations as shown in Figures \ref{fig:2-cycles}, \ref{fig:3-cycles}, \ref{fig:npc+}. Let $H(n)$, $n\ge 2$, be the kernel of the  homomorphism $G({\mathcal T},x^n)\to \mathbb Z_n$ that sends every generator to $1$. Then $H(n)$ is infinite and torsion-free. In particular all $G({\mathcal T, x^n})$, $n\ge 2$, are infinite.
\end{thm} 

{\bf Acknowledgements.} Work on this paper started as an REU project conducted at Boise State University in the Summer 2013, under the supervision of the second author. We gratefully acknowledge support from the Boise State University and the National Science Foundation. We also thank Stephan Rosebrock for many helpful discussions.

\section{Combinatorial Homotopy Theory and Curvature}
In this section we collect some facts on 2-dimensional combinatorial homotopy theory. Standard reference for diagrammatic methods in group theory and homotopy theory are the book \cite{MetzlerBook}, in particular Chapter V, written by Bogley and Pride \cite{BogleyPride}, and Gersten \cite{Gersten}. \\

Let $P=\langle x_1,...,x_k\ |\  r_1=1,...,r_m=1 \rangle$ be a presentation of a group $G$, and let $K(P)$ be the standard 2-complex obtained from that presentation.  A {\em surface diagram $S$ over $K(P)$} is a compact orientable surface $S$, with or without boundary, with a cell structure whose edges are oriented and labeled with the generators in $P$. Each 2-cell is marked by a base point positioned at a corner. We require that the word we obtain when reading around the boundary of a 2-cell in clockwise direction, starting at the base point, is either a relator $r_i$, or the inverse of a relator $r_i^{-1}$ of $P$. In the first case we mark the center of the cell with a plus sign, in the second case with a minus sign.\\

A map $f: L\to K$ between PLCW-complexes $L$ and $K$ is called {\em combinatorial} if it maps open cells homeomorphically to open cells.  A surface diagram $S$ over $K(P)$ gives rise to a combinatorial map $f: S\to K(P)$. The map $f$ sends all vertices in $S$ to the unique vertex in $K(P)$. An edge in $S$ labeled by a generator $x$ is send to the edge in $K(P)$ that corresponds to $x$ in an orientation preserving manner. A 2-cell in $S$ with boundary word $r^{\epsilon}$, $\epsilon =\pm 1$, is send to the 2-cell in $K(P)$ with that boundary word $r$, orientation preserving in case $\epsilon=1$ and orientation reversing in case $\epsilon=-1$. We can use {\em spherical diagrams} (surface diagrams where the underlying surface is the 2-sphere) to attach 3-cells to $K(P)$ by using the associated combinatorial maps as attaching maps.\\

\begin{figure}[here] 
   \centering
  
  \begin{tikzpicture}
\begin{scope}[decoration={markings, mark=at position 0.5 with {\arrow{>}}}] 
\draw [postaction={decorate}] (0,0) to (0,1);
\draw [postaction={decorate}] (0,1) to (1,1.8);
\draw [postaction={decorate}] (1,1.8) to (2,1);
\draw [postaction={decorate}] (2,1) to (2,0);
\draw [postaction={decorate}] (2,0) to (1,-0.8);
\draw [postaction={decorate}] (1,-0.8) to (0,0);

\draw [postaction={decorate}] (2,0) to (3,-0.8);
\draw [postaction={decorate}] (3,-0.8) to (4,0);
\draw [postaction={decorate}] (4,0) to (4,1);
\draw [postaction={decorate}] (4,1) to (3,1.8);
\draw [postaction={decorate}] (3,1.8) to (2,1);

\node [right] at (0,0.5) {$r$};
\node [below] at (0.5,1.3) {$r$};
\node [below] at (1.5,1.3) {$r$};
\node [left] at (2,0.5) {$r$};
\node [above] at (1.5,-0.3) {$r$};
\node [above] at (0.5,-0.3) {$r$};

\node [below] at (2.5,1.3) {$r$};
\node [below] at (3.5,1.3) {$r$};
\node [left] at (4,0.5) {$r$};
\node [above] at (3.5,-0.3) {$r$};
\node [above] at (2.5,-0.3) {$r$};
\node [right] at (2,0.5) {$r$};

\node [ ] at (1,0.5) {$+$};
\node [ ] at (3,0.5) {$-$};

\node [left] at (2,0.9) {$*$};
\node [right] at (2,0.1) {$*$};

\end{scope}
\end{tikzpicture}
 \caption{The spherical power diagram $S_{r^n}$ in case $n=6$. The two hexagons fold into a 2-sphere. The base points do not match and are one edge length apart.} 
   \label{fig:powerdiagram}
\end{figure}
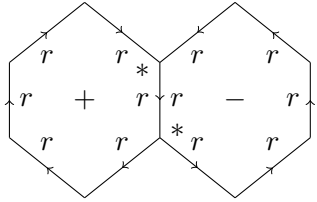

The second homotopy module $\pi_2(K(P))$ can be described entirely in terms of spherical diagrams. This is discussed in detail in \cite{BogleyPride}. However, Bogley and Pride prefer to work with {\em spherical pictures}, which are obtain from spherical diagrams by dualizing. If the presentation $P$ contains a proper power relation $r^n=1$, then $\pi_2(K(P)\ne 0$. A nontrivial $\pi_2$-element is shown in Figure \ref{fig:powerdiagram}. We call such $\pi_2$-elements {\em spherical power diagrams} and denote them by $S_{r^n}$. These spheres are referred to as primitive dipoles in \cite{BogleyPride}. \\ 

Let $K=K(P)$ be a standard 2-complex and $\bf A$ be a set of spherical diagrams over $K$. Given a spherical diagram $S$ over $K$ we can perform certain moves that transform $S$ into a new spherical diagram:

\begin{itemize}
\item Cancellation move: Suppose $S$ contains a vertex $v$ that is contained in the 2-cells $d$ and $d'$ and reading around $d$ in clockwise direction, starting at $v$, gives the same word $w$ as reading around $d'$ in anti clockwise direction. Then we can remove $d\cup d'$ from $S$ and close the resulting hole in $S$ by identifying along $w$. 
\item $\bf A$-replace move: Suppose that $S$ contains a disc diagram $D$ that is part of a spherical diagram $S'= D\cup D'$, where $S'$ or $-S'$, the mirror image of $S'$, is contained in $\bf A$. Then we can replace $D$ by $D'$ in $S$.
\end{itemize}

For a proof of the next result (in the language of pictures) see \cite{BogleyPride}, Section 1.2, page 162.

\begin{thm}\label{chfact} Let $K$ be a standard 2-complex constructed from a presentation $P$ of a group $G$ and let $\bf A$ be a set of spherical diagrams over $K$. Suppose every spherical diagram $S$ over $K$ can be transformed into the empty diagram using cancellation and $\bf A$-replace moves. Then the spherical diagrams in $\bf A$, together with the spherical power diagrams that come from the proper power relations in $P$, generate $\pi_2(K)$ as a $\mathbb ZG$-module. In particular, if we attach 3-cells to $K$ using the spherical diagrams in $\bf A$ and the spherical power diagrams as attaching maps, then the resulting 3-complex has trivial second homotopy module.
\end{thm}

We will also need the following result. See \cite{BogleyPride}, last paragraph on page 169. See also Huebschmann \cite{Huebschmann}.

\begin{thm}\label{maxfinitesubgroups} Let $K$ be a standard 2-complex constructed from a presentation $P$ of a group $G$. Suppose that $\pi_2(K)$ is generated by the spherical power diagrams that come from the proper power relations in $P$. If $r^n=1$ is a proper power relation in  $P$, then $\langle r \rangle$ is a subgroup of order $n$. Furthermore, the maximal finite subgroups of $G$ are cyclic and conjugate to the subgroups of the form $\langle r \rangle$, where $r^n=1$ is a proper power relation in $P$.
\end{thm}

We conclude this section with some facts concerning combinatorial curvature. Let $S$ be a compact orientable surface, possibly with boundary, with a cell structure. Assign real numbers $\alpha$ to the corners of the 2-cells in $S$. We think of these numbers as interior angles. The {\em curvature} at an interior vertex $v$ in $S$ is defined as $$\kappa(v)=2-\sum_i \alpha(v)_i,$$ where the $\alpha(v)_i$ range over the angles at $v$. If $v$ is a boundary vertex then we define $$\kappa(v)=1-\sum_i \alpha(v)_i.$$ The {\em curvature} of a 2-cell is defined as $$\kappa(d)=(\sum_j \alpha(d)_j)-(n-2)$$ where $n$ is the number of edges in the boundary of $d$ and the $\alpha(d)_j$ are the interior angles of $d$. The combinatorial Gauss-Bonnet theorem states that the total curvature is twice the Euler characteristic of the surface: $$\kappa(S)=\sum_v \kappa(v)+\sum_d\kappa(d)=2\chi(S).$$

\section{A 3-complex constructed from $\mathcal T$}

Let $\mathcal T$ be an injective labeled oriented tree that does not contain edge combinations as shown in Figures \ref{fig:2-cycles}, \ref{fig:3-cycles}, \ref{fig:npc+}. Let $$\bar P({\mathcal T},n)=\langle a, b, c,...\ |\ a^n=1, b^n=1, c^n=1,..., r_e=1, \{ e\in \mbox{edges of $\mathcal T$} \},  \rangle.$$ Let $\bar K$ be the standard 2-complex associated with $\bar P({\mathcal T},n)$. We refer to the 2-cells of $\bar K$ that come from the relators $r_e$ as {\em squares}, and to the 2-cells $d_{a^n}$, $d_{b^n}$, $d_{c^n}$,..., that come from the relators $a^n$, $b^n$, $c^n$,..., as {\em power 2-cells}. If $S$ is a surface diagram over $\bar K$, then a 2-cell in $S$ with boundary word $r_e^{\pm 1}$, is referred to as a {\em square in $S$}, and a 2-cell in $S$ with boundary word a power $a^{\pm n}$ is referred to as a {\em power 2-cell in $S$}.

\begin{figure}[htbp] 
   \centering
  
  \begin{tikzpicture}
\begin{scope}[decoration={markings, mark=at position 0.5 with {\arrow{>}}}] 
\draw [postaction={decorate}] (1,0) to (0.309, 0.951);
\draw [postaction={decorate}](0.309, 0.951) to (-0.808, 0.587);
\draw [postaction={decorate}](-0.808, 0.587) to (-0.809, -0.587);
\draw [postaction={decorate}](-0.809, -0.587) to (0.309, -0.951);
\draw [postaction={decorate}](0.309, -0.951) to (1,0);
\draw [postaction={decorate}] (2,0) to (2*0.309, 2*0.951);
\draw [postaction={decorate}](2*0.309, 2*0.951) to (-2*0.808, 2*0.587);
\draw [postaction={decorate}](2*-0.808, 2*0.587) to (2*-0.809, 2*-0.587);
\draw [postaction={decorate}](2*-0.809, 2*-0.587) to (2*0.309, 2*-0.951);
\draw [postaction={decorate}](2*0.309, 2*-0.951) to (2,0);
\draw [postaction={decorate}](1,0) to (2,0);
\draw [postaction={decorate}] (0.309, 0.951) to (2*0.309, 2*0.951);
\draw [postaction={decorate}]  (-0.808, 0.587) to (2*-0.808, 2*0.587);
\draw [postaction={decorate}] (-0.809, -0.587) to (2*-0.809, 2*-0.587);
\draw [postaction={decorate}] (0.309, -0.951) to (2*0.309, 2*-0.951);
\end{scope}
\node [above] at (1.5,0) {$b$};
\node [right] at  (1.5*0.309, 1.5*0.951) {$b$};
\node [right] at (1.5*-0.808, 1.5*0.587) {$b$};
\node [above] at (1.5*-0.809, 1.5*-0.587) {$b$};
\node [left] at  (1.5*0.309, 1.5*-0.951) {$b$};

\node [right] at (0.9*cos 36, 0.9*sin 36) {$a$};
\node [below] at (2*cos 36, 2*sin 36) {$c$};

\node [above] at (0.9*cos 108, 0.9*sin 108) {$a$};
\node [above] at (0.9*2*cos 108, 0.9*2*sin 108) {$c$};

\node [left] at (-0.9, 0) {$a$};
\node [left] at (-2*0.9, 0) {$c$};

\node [below] at (0.9*cos 252, 0.9*sin 252) {$a$};
\node [below] at (0.9*2*cos 252, 0.9*2*sin 252) {$c$};

\node [below] at (0.9*cos 324, 0.9*sin 324) {$a$};
\node [below] at (0.9*2*cos 324, 0.9*2*sin 324) {$c$};

\end{tikzpicture}
  
   \caption{The spherical diagram $S_e$ in case $n=5$. The top of the sphere is a 2-cell with boundary $c^n$. Base points are omitted.}
   \label{fig:S_e}
\end{figure}
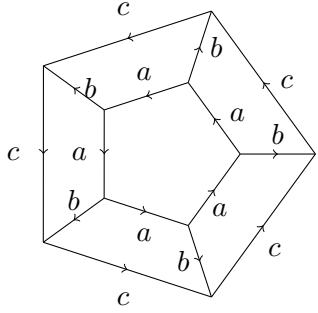

For every edge $e=(a|b|c)$ of $\mathcal T$ we have a spherical diagram $S_e$ with bottom a power 2-cell with boundary $a^n$, top a power 2-cell with boundary $c^n$, and sides a gallery of squares with boundary word $r_e$. See Figure \ref{fig:S_e}, and also Figure \ref{fig:3cellcollapse}. 

Let ${\bf A}=\{ S_e\ |\ e\in \mbox{edges of $\mathcal T$}\}$.  Let $L$ be the 3-complex obtained from $\bar K$ by attaching 3-cells using the spherical diagrams in $\bf A$ and spherical power diagrams that come from the proper power relations $a^n=1$, $b^n=1$, $c^n=1$,...., as attaching maps. We denote by $B_e$ the 3-cell in $L$ with attaching map $S_e$, and by $B_{a^n}$ the 3-cell in $L$ with attaching map the spherical power diagram that contains 2-cells with boundary $a^{\pm n}$.

\begin{thm}\label{thm:pi2} $\pi_2(L)=0$.
\end{thm}

Proof. In the light of Theorem \ref{chfact} it suffices to show that a spherical diagram $S$ over $\bar K$ can be transformed into the empty diagram using cancellation and $\bf A$-replace moves. Suppose that $S$ is a non-empty spherical diagram over $\bar K$ whose number of 2-cells can not be reduced using cancellation or $\bf A$-replace moves. 

\begin{lemma}\label{lemma:valency} If $v$ is a vertex in $S$ then $v$ has valency greater or equal to three.
\end{lemma}

Proof. The valency at $v$ can not be one because the relators in $\bar P({\mathcal T},n)$ are cyclically reduced. Suppose $v$ has valency two. Then $v$ is contained in exactly two 2-cells $d_1$ and $d_2$. These can not both be squares because that would imply that $\mathcal T$ contains an edge combination as shown in Figure \ref{fig:2-cycles}. If one of the 2-cells, say $d_1$, is a power 2-cell, then $d_2$ can not be a square because squares do not contain consecutive edges in their boundary with the same generator label. Thus both $d_1$ and $d_2$ are power 2-cells. But this implies that both $d_1$ and $d_2$ have the same boundary word and we can perform a cancellation. This contradicts the nature of the spherical diagram $S$. \qed \\

We assign to the corners of squares in $S$ the angle $1/2$ and to the corners of power cells in $S$ the angles $(n-2)/n$. Note that all 2-cells $d$ in $S$ are flat, that is $\kappa(d)=0$. Let $V$ be the set of vertices in $S$ and let $V_p$ be the subset consisting of vertices that are contained on the boundary of a power cell in $S$. If $v$ is not in $V_p$, then $v$ is surrounded by at least four squares because the complex $\bar K$ is a non-positively curved squared complex by hypothesis (we assumed that configurations shown in Figures \ref{fig:2-cycles}, \ref{fig:3-cycles} do not occur in $\mathcal T$). This implies $\kappa(v)\le 0$. It follows that $$4=\kappa(S)\le \sum_{v\in V_p} \kappa(v).$$

Let $v\in V_p$. Denote by $k(v)$ the number of power 2-cells in $S$ that contain $v$ in its boundary. Note that when traveling around the boundary of a power 2-cell in $S$ we encounter every vertex exactly once, because otherwise we would have $x^k=1$, where $x$ is some generator coming from a vertex in $\mathcal T$ and $k<n$. However, the fundamental group of $\bar K$ abelianizes to $\mathbb Z_n$, and $x$ maps to a generator of $\mathbb Z_n$. Thus $x$ has order $n$ in $\pi_1(\bar K)$. 

Define $$\tilde \kappa(v)=\kappa(v)/k(v).$$ Note that $$\sum_{v\in V_p} \kappa(v)=\sum_{d\in D_p}\sum_{v\in \partial d} \tilde \kappa(v),$$ where $D_p$ is the set of power 2-cells in $S$. In the following we will show that $\sum_{v\in \partial d} \tilde \kappa(v)\le 0$ for $d\in D_p$, which implies $4\le 0$, a contradiction. Consequently, a diagram $S$ with the stated properties does not exist.\\

\begin{lemma}\label{lemma:4squares} If $v$ is a vertex in $S$ of valency greater than three, then there have to be at least four squares among the 2-cells grouped around $v$ .
\end{lemma}

Proof. We have to have squares at $v$ because power 2-cells in $S$ do not share edges. There can not be a single square at $v$ because consecutive edges in squares carry different letters. Assume that there are exactly two squares at $v$. Since the valency at $v$ is greater than three, there have to be two power 2-cells at $v$, one with boundary $a^n$, and the other with boundary $b^n$. The situation is shown in Figure \ref{fig:2squares1}, orientations are omitted.   

\begin{figure}[here] 
   \centering

  \begin{tikzpicture}

\fill (0,0) circle (2pt); 
\draw (0,0) to (cos 45, sin 45);
\draw (0,0) to (cos 135, sin 135);
\draw (0,0) to (cos 225, sin 225);
\draw (0,0) to (cos 315,sin 315);
\draw (cos 45, sin 45) to (sqrt 2*cos 90, sqrt 2 * sin 90);
\draw (cos 135, sin 135) to (sqrt 2*cos 90, sqrt 2 * sin 90);
\draw (cos 225, sin 225) to (sqrt 2*cos -90, sqrt 2 * sin -90);
\draw (cos 315, sin 315) to (sqrt 2*cos -90, sqrt 2 * sin -90);

\node [ ] at (0.5*cos 60, 0.5*sin 60){$b$};
\node [ ] at (0.5*cos 120, 0.5*sin 120){$a$};
\node [ ] at (0.5*cos 300, 0.5*sin 300){$b$};
\node [ ] at (0.5*cos 240, 0.5*sin 240){$a$};

\end{tikzpicture}
 \caption{Two squares at $v$. Between edges with the same letter is a power 2-cell.} 
   \label{fig:2squares1}
\end{figure}
Let $s_1$ be the top square and $s_2$ be the bottom square, and let $e_i$ be the edges in $\mathcal T$ associated with $s_i$, $i=1,2$. Note that either $a$ or $b$ occurs twice in the boundary of $s_1$, so either $a$ or $b$ is the label on $e_1$. If the two edges $e_1$ and $e_2$ are different, then $a$ and $b$ occur as labels on these two edges, because we assumed $\mathcal T$ to be injective. But this implies that $\mathcal T$ contains a edge combination as shown in Figure \ref{fig:2-cycles}. Thus $e_1=e_2$ and we assume without loss of generality that the label on $e_1$ is $a$. We are in the situation shown in Figure \ref{fig:2squares}.
\begin{figure}[here] 
   \centering
    \begin{tikzpicture}
\fill (0,0) circle (2pt); 
\draw (0,0) to (cos 45, sin 45);
\draw (0,0) to (cos 135, sin 135);
\draw (0,0) to (cos 225, sin 225);
\draw (0,0) to (cos 315,sin 315);
\draw (cos 45, sin 45) to (sqrt 2*cos 90, sqrt 2 * sin 90);
\draw (cos 135, sin 135) to (sqrt 2*cos 90, sqrt 2 * sin 90);
\draw (cos 225, sin 225) to (sqrt 2*cos -90, sqrt 2 * sin -90);
\draw (cos 315, sin 315) to (sqrt 2*cos -90, sqrt 2 * sin -90);
\node [ ] at (0.5*cos 60, 0.5*sin 60){$b$};
\node [ ] at (0.5*cos 120, 0.5*sin 120){$a$};
\node [ ] at (0.5*cos 300, 0.5*sin 300){$b$};
\node [ ] at (0.5*cos 240, 0.5*sin 240){$a$};
\node [ ] at (0.5*cos 120 + cos 45, 0.5*sin 120 + sin 45){$a$};
\node [ ] at (0.5*cos 240 + cos -45, 0.5*sin 240 + sin -45){$a$};
\node [ ] at (0.5*cos 60 + cos 135, 0.5*sin 60 + sin 135){$c$};
\node [ ] at (0.5*cos 300 + cos 225 , 0.5*sin 300 + sin 225){$c$};
\end{tikzpicture}
 \caption{Two squares at $v$. Between edges with the same letter is a power 2-cell.} 
   \label{fig:2squares}
\end{figure}
But this leads to contradicting orientations on the edge $e_1$ that gives rise to the squares. Let us discuss one case in detail. 
Assume that the edges labeled $a$ and $b$ in $s_1$ that contain $v$ point towards $v$. Then the edges with these labels in $s_2$ that contain $v$ point away from $v$, because the $a$-edges that contain $v$ are part of a power 2-cell and so are the $b$-edges. But then the edge associated with $s_1$ is $(c|a|b)$ and the edge associated with $s_2$ is $(b|a|c)$. These two edges have the same vertices and the same label, but have opposite orientation and hence are not the same. The other cases lead to the same contradiction. 

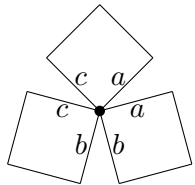
\begin{figure}[here] 
   \centering
  
  \begin{tikzpicture}
  \fill (0,0) circle (2pt);
 
\draw (0,0) to (cos 45, sin 45);
\draw (0,0) to (cos 135, sin 135);
\draw (0,0) to (cos 165, sin 165);
\draw (0,0) to (cos 255,sin 255);
\draw (0,0) to (cos 285,sin 285);
\draw (0,0) to (cos 15,sin 15);
\draw (cos 45, sin 45) to (sqrt 2*cos 90, sqrt 2*sin 90);
\draw (cos 135, sin 135) to (sqrt 2*cos 90, sqrt 2*sin 90);
\draw (cos 165, sin 165) to (sqrt 2*cos 210, sqrt 2*sin 210);
\draw (cos 255, sin 255) to (sqrt 2*cos 210, sqrt 2*sin 210);
\draw (cos 285,sin 285) to (sqrt 2*cos 330, sqrt 2*sin 330);
\draw (cos 15,sin 15) to (sqrt 2*cos 330, sqrt 2*sin 330);
\node [] at (0.5*cos 60, 0.5*sin 60) {$a$};
\node [] at (0.5*cos 0, 0.5*sin 0) {$a$};
\node [] at (0.5*cos 300, 0.5*sin 300) {$b$};
\node [] at (0.5*cos 240, 0.5*sin 240) {$b$};
\node [] at (0.5*cos 180, 0.5*sin 180) {$c$};
\node [] at (0.5*cos 120, 0.5*sin 120) {$c$};
\end{tikzpicture}
 \caption{Three squares around $v$. Between edges with the same letter there may or may not be a power 2-cell. If not the edges are to be identified.} 
   \label{fig:3squares}
\end{figure}

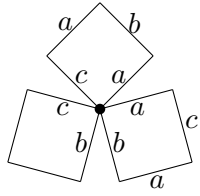
\begin{figure}[here] 
   \centering
  
  \begin{tikzpicture}
  \fill (0,0) circle (2pt);
 \draw (0,0) to (cos 45, sin 45);
\draw (0,0) to (cos 135, sin 135);
\draw (0,0) to (cos 165, sin 165);
\draw (0,0) to (cos 255,sin 255);
\draw (0,0) to (cos 285,sin 285);
\draw (0,0) to (cos 15,sin 15);
\draw (cos 45, sin 45) to (sqrt 2*cos 90, sqrt 2*sin 90);
\draw (cos 135, sin 135) to (sqrt 2*cos 90, sqrt 2*sin 90);
\draw (cos 165, sin 165) to (sqrt 2*cos 210, sqrt 2*sin 210);
\draw (cos 255, sin 255) to (sqrt 2*cos 210, sqrt 2*sin 210);
\draw (cos 285,sin 285) to (sqrt 2*cos 330, sqrt 2*sin 330);
\draw (cos 15,sin 15) to (sqrt 2*cos 330, sqrt 2*sin 330);

\node [] at (0.5*cos 60, 0.5*sin 60) {$a$};
\node [] at (0.5*cos 0, 0.5*sin 0) {$a$};
\node [] at (0.5*cos 300, 0.5*sin 300) {$b$};
\node [] at (0.5*cos 240, 0.5*sin 240) {$b$};
\node [] at (0.5*cos 180, 0.5*sin 180) {$c$};
\node [] at (0.5*cos 120, 0.5*sin 120) {$c$};

\node [] at (0.5*cos 60 + cos 135, 0.5*sin 60 + sin 135) {$a$};
\node [] at (0.5*cos 0 + cos 285, 0.5*sin 0 + sin 285) {$a$};
\node [] at (0.5*cos 300 + cos 15, 0.5*sin 300 + sin 15) {$c$};
\node [] at (0.5*cos 120 + cos 45, 0.5*sin 120 + sin 45) {$b$};

\end{tikzpicture}
\caption{Three squares around $v$ with two squares arising from the same edge in $\mathcal T$.} 
\label{fig:3squares+}
\end{figure}

Let us assume next that there are exactly three squares at $v$. There can not be more than one power 2-cell between two consecutive squares because power 2-cell in $S$ do not share edges. There has to be at least one power 2-cell at $v$, because we assume the valency at $v$ is greater than three. The situation is shown in Figure \ref{fig:3squares}, orientations are omitted. 
Note that the letters $a$, $b$, and $c$ have to be all different, because there is no square that has the same letter on consecutive edges.  Let $s_1$ be the square with sides $a$ and $b$, $s_2$ be the square with sides $b$ and $c$, and $s_3$ be the square with sides $a$ and $c$. Let $e_i$ be the edge in $\mathcal T$ associated with $s_i$, $i=1,2,3$. If all three edges are different, then, since we assumed ${\mathcal T}$ to be injective,  $a$, $b$, and $c$ occur as labels on these three edges. But this implies that $\mathcal T$ contains a edge combination as shown in Figure \ref{fig:3-cycles}. Thus two squares come from the same edge. Without loss of generality we assume that $e_1=e_3$, and $e_1$ has label $a$. We are in the situation shown in Figure \ref{fig:3squares+}. 
But this implies that $\mathcal T$ contains an edge combination as shown in Figure \ref{fig:npc+}. Thus there have to be at least four squares at $v$. \qed

\begin{lemma}\label{lem:curvature} Let $v$ be a vertex in $S$ that is contained in a power 2-cell in $S$.
\begin{itemize}
\item If the valency at $v$ is three then $\tilde\kappa(v)=\frac{2}{n}$.
\item If the valency at $v$ is greater than three then $\tilde\kappa(v)\le -1+\frac{2}{n}\le 0$.
\end{itemize}
\end{lemma}

Proof. If $v$ has valency three then $v$ is contained in exactly three 2-cells $d_1$, $d_2$, $d_3$ in $S$. Not all of these can be squared because that would imply that $\mathcal T$ contains an edge combination as shown in Figure \ref{fig:3-cycles}. Since any pair among these share an edge and we assumed that no cancellation can take place in $S$, there can only be one power 2-cell among  the $d_i$, $i=1,2,3$. Thus we may assume that $d_1$ is a power 2-cell and $d_2$, $d_3$ are squares. This gives $$\tilde\kappa(v)=2-(2\cdot\frac{1}{2}+\frac{n-2}{n})=\frac{2}{n}.$$\\
 
Next assume that $v$ has valency greater than 3. By Lemma \ref{lemma:4squares} there are at least four squares grouped around $v$. Thus we have $$\kappa(v)\le 2-(4\cdot\frac{1}{2}+k(v)\cdot\frac{n-2}{n}))=-k(v)\cdot\frac{n-2}{n}.$$ Since $v$ is contained in a power 2-cell $k(v)\ge 1$ and we have $$\tilde\kappa(v)\le -1+\frac{2}{n}\le 0.$$ \qed









\begin{lemma}\label{lemma:val3} Let $v$ be a vertex in $S$ of valency three and assume that $v$ is in the boundary of a power 2-cell. Then the local configuration we have in $S$ at $v$ also occurs in some spherical diagram $S_e$ or $-S_e$, the mirror image of $S_e$.
\end{lemma}

\begin{figure}[htbp] 
   \centering
    
   \begin{tikzpicture}
\begin{scope}[decoration={markings, mark=at position 0.5 with {\arrow{>}}}] 
\draw [postaction={decorate}] (1,0) to (0.309, 0.951);
\draw [postaction={decorate}](0.309, 0.951) to (-0.808, 0.587);
\draw [postaction={decorate}](-0.808, 0.587) to (-0.809, -0.587);
\draw [postaction={decorate}](-0.809, -0.587) to (0.309, -0.951);
\draw [postaction={decorate}](0.309, -0.951) to (1,0);
\draw [postaction={decorate}] (2,0) to (2*0.309, 2*0.951);
\draw [postaction={decorate}](2*0.309, 2*0.951) to (-2*0.808, 2*0.587);
\draw [postaction={decorate}](1,0) to (2,0);
\draw [postaction={decorate}] (0.309, 0.951) to (2*0.309, 2*0.951);
\draw [postaction={decorate}]  (-0.808, 0.587) to (2*-0.808, 2*0.587);
\end{scope}

\node [right] at  (1.5*0.309, 1.5*0.951) {$b$};
\node [right] at (0.9*cos 36, 0.9*sin 36) {$a$};
\node [above] at (0.9*cos 108, 0.9*sin 108) {$a$};
\node [left] at (-0.9, 0) {$a$};
\node [below] at (0.9*cos 252, 0.9*sin 252) {$a$};
\node [below] at (0.9*cos 324, 0.9*sin 324) {$a$};
\end{tikzpicture}
 \caption{A vertex of valency 3 in $S$.}
   \label{fig:val3}
\end{figure}
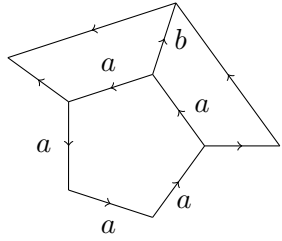

Proof. If $v$ were in the boundary of two power 2-cells, then the valency at $v$ would be at least four because power 2-cells in $S$ can not share edges. So $v$ is surrounded by one power 2-cell and two squares and we are in the situation shown in Figure \ref{fig:val3}. We assume the orientation on the power 2-cell is counter-clockwise and the edge labeled $b$ points away from the vertex. The other cases can be argued in exactly the same fashion.  
Note that if the two squares come from different edges $e_1$ and $e_2$ in $\mathcal T$, then one has label $a$ and the other has label $b$ (injectivity of $\mathcal T$), but that implies that the two edges $e_1$ and $e_2$ in $\mathcal T$ are as shown in Figure \ref{fig:2-cycles}. Thus the two squares come from a single edge $e$. Suppose that $a$ is the edge label on $e$ in $\mathcal T$. This yields to a contradiction concerning the orientation of $e$. One square would imply that $e=(b|a|c)$ the other would imply that $e=(c|a|b)$. Thus $b$ is the label on $e$ and local configuration as shown in Figure \ref{fig:val3} is indeed also present in the spherical diagram $S_e$ as shown in Figure \ref{fig:S_e}.\qed \\ 

We will continue with the proof of Theorem \ref{thm:pi2}. Consider a power cell $d$ in $S$. Assume first that $n$ is even. If we have $n/2$ or more vertices of valency 3 in the boundary of $d$ in succession, then, by Lemma \ref{lemma:val3},  $S$ contains more than half of some spherical diagram $\pm S_e$. In that case we could use an $\bf A$-replace move to reduce the number of 2-cells in $S$, which we assumed not to be possible. Thus the boundary of $d$ can contain at most $(n/2)-1$ vertices of valency 3 in succession. Similarly, if $n$ is odd and we have $(n-1)/2$ or more vertices of valency 3 in the boundary of $d$ is succession, we could use an $\bf A$-replace move to reduce the number of 2-cells in $S$, which we assumed not to be possible. Thus the boundary of $d$ can contain at most $(n-1)/2-1$ vertices of valency 3 in succession. So, no matter if $n$ is even or odd, the boundary of $d$ must contain at least two vertices of valency greater than three. By Lemma \ref{lem:curvature} we have $$\sum_{v\in \partial d}\tilde \kappa(v)\le(n-2)\cdot\frac{2}{n}+2\cdot (-1+\frac{2}{n}) = 0.$$ It follows that $$4=\kappa(S)\le 0,$$ a contradiction. We conclude that there does not exist a non-empty spherical diagram over $L^{(2)}$ whose number of 2-cells can not be reduced by cancellation and $\bf A$-replace moves. \qed

\begin{thm}\label{thm:notcyclic} $G=\pi_1(L)$ is not cyclic.
\end{thm}

Proof. Assume that $G$ is cyclic. Since we can map $G$ to $\mathbb Z_n$ by sending every generator that comes from a vertex of $\mathcal T$ to $1$, the group $G$ is cyclic of order $n$ and any vertex of $\mathcal T$ gives the same generator of $G$. Choose two vertices $g$, $h$ of $\mathcal T$ that are not contained in the set of vertices and labels of a single edge in $\mathcal T$. Let us make sure that such a pair of vertices exists. Choose a pair of edges $e_1$ and $e_2$ in $\mathcal T$ that share a vertex $b$. Let $a$ and $c$ be the other vertices of $e_1$ and $e_2$, respectively. Note that the label $x$ on $e_1$ is not $a$ or $b$ because $\mathcal T$ is compressed, and it is not $c$ because no edge combination as shown in Figure \ref{fig:npc+} occurs in $\mathcal T$. For the same reason the label $y$ on $e_2$ is not $a$, $b$, or $c$. Note that the labels $x$ and $y$ are distinct and can not occur elsewhere as edge labels because $\mathcal T$ is injective. If there is no edge with vertices $y$ and $a$, choose $g=a$, $h=y$. If there is such an edge, but there is no edge with vertices $c$ and $x$, choose $g=c$, $h=x$. If there is such an edge, choose $g=x$, $h=y$.
The vertices $g$ and $h$ are chosen in this way so that, if $g$ and $h$ are considered as generators of $G$, they do not both occur on the boundary of a single square $d_e$ of $L^{(2)}=\bar K$. Since $g=h$ in $G$, there exists a disc diagram $E$ over $L^{(2)}$ with boundary $gh^{-1}$. We assume that the number of 2-cells in $E$ can not be reduced using cancellation or $\bf A$-replace moves. We will show, using essentially the same reasoning as used in the proof of Theorem \ref{thm:pi2}, that $\kappa(E)\le 0$. A contradiction since $E$ is a disc and hence $\kappa(E)=2$. 

We use the same notation as in the proof of Theorem \ref{thm:pi2}. Note that we have exactly two vertices on the boundary of $E$.  Let $V$ for the set of vertices of $E$ and $V_p$ for the set of vertices that are contained in the boundary of a power 2-cell in $E$. Consider a vertex of $V-V_p$. If $v$ is an interior vertex then there are at least four squares grouped around $v$ because we assumed $\mathcal T$ does not contain edge combinations as shown in Figures \ref{fig:2-cycles}, \ref{fig:3-cycles}. Thus $\kappa(v)\le 0$. Assume next that $v$ is a boundary vertex. Since we assumed that $g$ and $h$ are not contained in the boundary of a single square, there are at least two squares grouped around $v$. Thus $\kappa(v)\le 0$. It follows that $$\kappa(E)\le \sum_{v\in V_p}\kappa(v)=\sum_{d\in D_p}\sum_{v\in \partial d}\tilde\kappa(v),$$ where $D_p$ denotes the set of power 2-cells in $E$. We will now proceed as in the proof of Theorem \ref{thm:pi2} to show that $$\sum_{v\in \partial d}\tilde\kappa(v)\le 0$$ for every power 2-cell $d$ in $E$. Let $d$ be a power 2-cell in $E$. If all vertices of $d$ are interior vertices then the sum in question is indeed less or equal to $0$ by exactly the same reasoning as used in the proof of Theorem \ref{thm:pi2}. So suppose $d$ does contain a vertex $v$ on the boundary of $E$. Again, since we assumed that $g$ and $h$ are not contained in the boundary of a single square, there are at least two squares grouped around $v$. Thus $$\kappa(v)\le 1-(2\cdot \frac{1}{2}+k(v)\cdot\frac{(n-2)}{n})$$ and hence $$\tilde\kappa(v)\le -1+\frac{2}{n}.$$ Note that this is the same inequality that we have for interior vertices of valency greater that three in $E$ according to Lemma \ref{lem:curvature}. Not all vertices on the boundary of $d$ that are different from the boundary vertex $v$ can be interior vertices of valency three, because in that case we could perform an $\bf A$-replace move that would lower the number of 2-cells in $E$. So $\partial d$ contains a vertex $w$ distinct from $v$ that is either interior and has valency greater than three, or is the other boundary vertex of $E$. Thus $$\sum_{v\in \partial d}\tilde \kappa(v)\le(n-2)\cdot\frac{2}{n}+2\cdot (-1+\frac{2}{n}) = 0.$$
We have reached the desired contradiction. The assumption that $G$ is cyclic is false. \qed 

\section{Collapsing the 3-complex $L$}

Let $\mathcal T$ be the labeled oriented tree that was used for the construction of the 3-complex $L$. The tree $\mathcal T$ can be collapsed to a single vertex. We start with and edge $e$ that contains a boundary vertex, i.e. a vertex of valency one in $\mathcal T$. Let us assume $e=(a|b|c)$ and $c$ is the boundary vertex. The vertex $c$ is a free face of $e$, and we can collapse $e$ by pushing in $c$. This removes the edge $e$ and the vertex $c$ from $\mathcal T$ and leaves us with a labeled oriented tree ${\mathcal T}'$. The vertex $c$ might show up as an edge label in ${\mathcal T}'$ despite the fact that ${\mathcal T}'$ does not contain the vertex $c$, but this will not concern us. We now choose an edge in ${\mathcal T}'$ that contains a boundary vertex and collapse it in the described fashion. We continue until $\mathcal T$ is collapsed to a single vertex. It is not difficult to see that one can choose a vertex $z$ of $\mathcal T$ and then devise a collapsing strategy that collapses $\mathcal T$ to the chosen vertex $z$.

\begin{figure}[htbp] 
\centering
\begin{tikzpicture}
\fill (0,4.7) circle (2pt);
\node[below] at (0,4.7) {$a$};
\fill (4,4.7) circle (2pt);
\node[below] at (4,4.7) {$c$};

\begin{scope}[decoration={markings, mark=at position 0.5 with {\arrow{>}}}]
\draw [postaction={decorate}] (0,4.7) --(4,4.7) node[midway, above]{$b$};
\end{scope}

\begin{scope}[decoration={markings, mark=at position 0.5 with {\arrow{>}}}]
\draw [postaction={decorate}] (0,0) --(0.5,0.85) node[midway, right]{$a$};
\draw [postaction={decorate}] (0.5,0.85) --(0.5,1.85) node[midway, right]{$a$};
\draw [postaction={decorate}] (0.5,1.85) --(0,1.85+0.85) node[midway, right]{$a$};
\end{scope}
\begin{scope}[decoration={markings, mark=at position 0.5 with {\arrow{<}}}]
\draw [postaction={decorate}] (0,0) -- (-0.5,0.85) node[midway, right]{$a$};
\draw [postaction={decorate}] (-0.5,0.85) -- (-0.5,1.85) node[midway, right]{$a$};
\draw [postaction={decorate}] (-0.5,1.85) -- (0,1.85+0.85) node[midway, right]{$a$};
\end{scope}

\begin{scope}[decoration={markings, mark=at position 0.5 with {\arrow{>}}}]
\draw [postaction={decorate}] (0,0) -- (4,0) node[midway, above]{$b$};
\draw [postaction={decorate}] (0.5,0.85) -- (4.5, 0.85) node[midway, above]{$b$};
\draw [postaction={decorate}] (0.5,1.85) -- (4.5, 1.85) node[midway, above]{$b$};
\draw [postaction={decorate}] (0,2.7) -- (4,2.7) node[midway, above]{$b$};
\end{scope}

\begin{scope}[decoration={markings, mark=at position 0.5 with {\arrow{>}}}]
\draw [postaction={decorate}] (4,0) -- (4.5,0.85) node[midway, right]{$c$};
\draw [postaction={decorate}] (4.5,0.85) -- (4.5,1.85) node[midway, right]{$c$};
\draw [postaction={decorate}] (4.5,1.85) -- (4,1.85+0.85) node[midway, right]{$c$};
\end{scope}
\begin{scope}[decoration={markings, mark=at position 0.5 with {\arrow{<}}}]
\draw [dashed, postaction={decorate}] (4,0) -- (3.5,0.85) node[midway, right]{$c$};
\draw [dashed, postaction={decorate}] (3.5,0.85) -- (3.5,1.85) node[midway, right]{$c$};
\draw [dashed, postaction={decorate}] (3.5,1.85) -- (4,1.85+0.85) node[midway, right]{$c$};
\end{scope}
\end{tikzpicture}
 \caption{The edge $e=(a|b|c)$ and the corresponding 3-cell $B_e$ of $L$. We can internally collapse $B_e$ by pushing in the 2-cell $d_{c^n}$ which is part of the boundary of $B_e$.}
   \label{fig:3cellcollapse}
\end{figure}
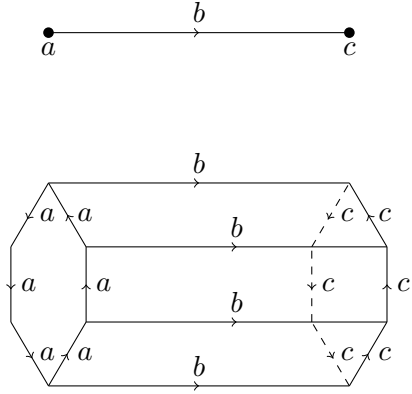
\begin{figure}[htbp] 
\centering
\begin{tikzpicture}[scale=0.8]

\node at (1.8*cos 30,1.8*sin 30){$*$};

\begin{scope}[decoration={markings, mark=at position 0.5 with {\arrow{<}}}]
\draw [postaction={decorate}] (2*cos -30,2*sin -30) --(2*cos 30,2*sin 30) node[midway, left]{$c$};
\end{scope}

\begin{scope}[decoration={markings, mark=at position 0.5 with {\arrow{<}}}]
\draw [postaction={decorate}] (2*cos 30,2*sin 30) --(2*cos 90,2*sin 90) node[midway, below]{$c$};
\end{scope}

\begin{scope}[decoration={markings, mark=at position 0.5 with {\arrow{<}}}]
\draw [postaction={decorate}] (2*cos 90,2*sin 90) --(2*cos 150,2*sin 150) node[midway, below]{$c$};
\end{scope}

\begin{scope}[decoration={markings, mark=at position 0.5 with {\arrow{<}}}]
\draw [postaction={decorate}] (2*cos 150,2*sin 150) --(2*cos 210,2*sin 210) node[midway, right]{$c$};
\end{scope}

\begin{scope}[decoration={markings, mark=at position 0.5 with {\arrow{<}}}]
\draw [postaction={decorate}] (2*cos 210,2*sin 210) --(2*cos 270,2*sin 270) node[midway, above]{$c$};
\end{scope}

\begin{scope}[decoration={markings, mark=at position 0.5 with {\arrow{<}}}]
\draw [postaction={decorate}] (2*cos 270,2*sin 270) --(2*cos 330,2*sin 330) node[midway, above]{$c$};
\end{scope}

\begin{scope}[xshift=3.4641cm]

\node at (1.8*cos 210,1.8*sin 210){$*$};

\begin{scope}[decoration={markings, mark=at position 0.5 with {\arrow{>}}}]
\draw [postaction={decorate}] (2*cos -30,2*sin -30) --(2*cos 30,2*sin 30) node[midway, left]{$c$};
\end{scope}

\begin{scope}[decoration={markings, mark=at position 0.5 with {\arrow{>}}}]
\draw [postaction={decorate}] (2*cos 30,2*sin 30) --(2*cos 90,2*sin 90) node[midway, below]{$c$};
\end{scope}

\begin{scope}[decoration={markings, mark=at position 0.5 with {\arrow{>}}}]
\draw [postaction={decorate}] (2*cos 90,2*sin 90) --(2*cos 150,2*sin 150) node[midway, below]{$c$};
\end{scope}

\begin{scope}[decoration={markings, mark=at position 0.5 with {\arrow{>}}}]
\draw [postaction={decorate}] (2*cos 210,2*sin 210) --(2*cos 270,2*sin 270) node[midway, above]{$c$};
\end{scope}

\begin{scope}[decoration={markings, mark=at position 0.5 with {\arrow{>}}}]
\draw [postaction={decorate}] (2*cos 270,2*sin 270) --(2*cos 330,2*sin 330) node[midway, above]{$c$};
\end{scope}

\end{scope}


\begin{scope}[xshift=7.5cm]

\node at (0.8*cos 30,0.8*sin 30){$*$};

\begin{scope}[decoration={markings, mark=at position 0.5 with {\arrow{<}}}]
\draw [postaction={decorate}] (2*cos -30,2*sin -30) --(2*cos 30,2*sin 30) node[midway, left]{$c$};
\end{scope}

\begin{scope}[decoration={markings, mark=at position 0.5 with {\arrow{<}}}]
\draw [postaction={decorate}] (2*cos 30,2*sin 30) --(2*cos 90,2*sin 90) node[midway, below]{$c$};
\end{scope}

\begin{scope}[decoration={markings, mark=at position 0.5 with {\arrow{<}}}]
\draw [postaction={decorate}] (2*cos 90,2*sin 90) --(2*cos 150,2*sin 150) node[midway, below]{$c$};
\end{scope}

\begin{scope}[decoration={markings, mark=at position 0.5 with {\arrow{<}}}]
\draw [postaction={decorate}] (2*cos 150,2*sin 150) --(2*cos 210,2*sin 210) node[midway, right]{$c$};
\end{scope}

\begin{scope}[decoration={markings, mark=at position 0.5 with {\arrow{<}}}]
\draw [postaction={decorate}] (2*cos 210,2*sin 210) --(2*cos 270,2*sin 270) node[midway, above]{$c$};
\end{scope}

\begin{scope}[decoration={markings, mark=at position 0.5 with {\arrow{<}}}]
\draw [postaction={decorate}] (2*cos 270,2*sin 270) --(2*cos 330,2*sin 330) node[midway, above]{$c$};
\end{scope}

\begin{scope}[decoration={markings, mark=at position 0.5 with {\arrow{<}}}]
\draw [postaction={decorate}] (cos -30,sin -30) --(cos 30,sin 30) node[midway, left]{$a$};
\end{scope}

\begin{scope}[decoration={markings, mark=at position 0.5 with {\arrow{<}}}]
\draw [postaction={decorate}] (cos 30,sin 30) --(cos 90,sin 90) node[midway, below]{$a$};
\end{scope}

\begin{scope}[decoration={markings, mark=at position 0.5 with {\arrow{<}}}]
\draw [postaction={decorate}] (cos 90,sin 90) --(cos 150,sin 150) node[midway, below]{$a$};
\end{scope}

\begin{scope}[decoration={markings, mark=at position 0.5 with {\arrow{<}}}]
\draw [postaction={decorate}] (cos 150,sin 150) --(cos 210,sin 210) node[midway, right]{$a$};
\end{scope}

\begin{scope}[decoration={markings, mark=at position 0.5 with {\arrow{<}}}]
\draw [postaction={decorate}] (cos 210,sin 210) --(cos 270,sin 270) node[midway, above]{$a$};
\end{scope}

\begin{scope}[decoration={markings, mark=at position 0.5 with {\arrow{<}}}]
\draw [postaction={decorate}] (cos 270,sin 270) --(cos 330,sin 330) node[midway, above]{$a$};
\end{scope}

\begin{scope}[decoration={markings, mark=at position 0.5 with {\arrow{>}}}]
\draw [postaction={decorate}] (cos -30,sin -30) --(2*cos -30,2*sin -30) node[midway, below]{$b$};
\end{scope}

\begin{scope}[decoration={markings, mark=at position 0.5 with {\arrow{>}}}]
\draw [postaction={decorate}] (cos 30,sin 30) --(2*cos 30,2*sin 30) node[midway, below]{$b$};
\end{scope}

\begin{scope}[decoration={markings, mark=at position 0.5 with {\arrow{>}}}]
\draw [postaction={decorate}] (cos 90,sin 90) --(2*cos 90,2*sin 90) node[midway, right]{$b$};
\end{scope}

\begin{scope}[decoration={markings, mark=at position 0.5 with {\arrow{>}}}]
\draw [postaction={decorate}] (cos 150,sin 150) --(2*cos 150,2*sin 150) node[midway, above]{$b$};
\end{scope}

\begin{scope}[decoration={markings, mark=at position 0.5 with {\arrow{>}}}]
\draw [postaction={decorate}] (cos 210,sin 210) --(2*cos 210,2*sin 210) node[midway, right]{$b$};
\end{scope}

\begin{scope}[decoration={markings, mark=at position 0.5 with {\arrow{>}}}]
\draw [postaction={decorate}] (cos 270,sin 270) --(2*cos 270,2*sin 270) node[midway, right]{$b$};
\end{scope}

\begin{scope}[xshift=3.4641cm]

\node at (0.8*cos 210,0.8*sin 210){$*$};

\begin{scope}[decoration={markings, mark=at position 0.5 with {\arrow{>}}}]
\draw [postaction={decorate}] (2*cos -30,2*sin -30) --(2*cos 30,2*sin 30) node[midway, left]{$c$};
\end{scope}

\begin{scope}[decoration={markings, mark=at position 0.5 with {\arrow{>}}}]
\draw [postaction={decorate}] (2*cos 30,2*sin 30) --(2*cos 90,2*sin 90) node[midway, below]{$c$};
\end{scope}

\begin{scope}[decoration={markings, mark=at position 0.5 with {\arrow{>}}}]
\draw [postaction={decorate}] (2*cos 90,2*sin 90) --(2*cos 150,2*sin 150) node[midway, below]{$c$};
\end{scope}


\begin{scope}[decoration={markings, mark=at position 0.5 with {\arrow{>}}}]
\draw [postaction={decorate}] (2*cos 210,2*sin 210) --(2*cos 270,2*sin 270) node[midway, above]{$c$};
\end{scope}

\begin{scope}[decoration={markings, mark=at position 0.5 with {\arrow{>}}}]
\draw [postaction={decorate}] (2*cos 270,2*sin 270) --(2*cos 330,2*sin 330) node[midway, above]{$c$};
\end{scope}

\begin{scope}[decoration={markings, mark=at position 0.5 with {\arrow{>}}}]
\draw [postaction={decorate}] (cos -30,sin -30) --(cos 30,sin 30) node[midway, left]{$a$};
\end{scope}

\begin{scope}[decoration={markings, mark=at position 0.5 with {\arrow{>}}}]
\draw [postaction={decorate}] (cos 30,sin 30) --(cos 90,sin 90) node[midway, below]{$a$};
\end{scope}

\begin{scope}[decoration={markings, mark=at position 0.5 with {\arrow{>}}}]
\draw [postaction={decorate}] (cos 90,sin 90) --(cos 150,sin 150) node[midway, below]{$a$};
\end{scope}

\begin{scope}[decoration={markings, mark=at position 0.5 with {\arrow{>}}}]
\draw [postaction={decorate}] (cos 150,sin 150) --(cos 210,sin 210) node[midway, right]{$a$};
\end{scope}

\begin{scope}[decoration={markings, mark=at position 0.5 with {\arrow{>}}}]
\draw [postaction={decorate}] (cos 210,sin 210) --(cos 270,sin 270) node[midway, above]{$a$};
\end{scope}

\begin{scope}[decoration={markings, mark=at position 0.5 with {\arrow{>}}}]
\draw [postaction={decorate}] (cos 270,sin 270) --(cos 330,sin 330) node[midway, above]{$a$};
\end{scope}

\begin{scope}[decoration={markings, mark=at position 0.5 with {\arrow{>}}}]
\draw [postaction={decorate}] (cos -30,sin -30) --(2*cos -30,2*sin -30) node[midway, below]{$b$};
\end{scope}

\begin{scope}[decoration={markings, mark=at position 0.5 with {\arrow{>}}}]
\draw [postaction={decorate}] (cos 30,sin 30) --(2*cos 30,2*sin 30) node[midway, below]{$b$};
\end{scope}

\begin{scope}[decoration={markings, mark=at position 0.5 with {\arrow{>}}}]
\draw [postaction={decorate}] (cos 90,sin 90) --(2*cos 90,2*sin 90) node[midway, right]{$b$};
\end{scope}

\begin{scope}[decoration={markings, mark=at position 0.5 with {\arrow{>}}}]
\draw [postaction={decorate}] (cos 150,sin 150) --(2*cos 150,2*sin 150) node[midway, above]{$b$};
\end{scope}

\begin{scope}[decoration={markings, mark=at position 0.5 with {\arrow{>}}}]
\draw [postaction={decorate}] (cos 210,sin 210) --(2*cos 210,2*sin 210) node[midway, right]{$b$};
\end{scope}

\begin{scope}[decoration={markings, mark=at position 0.5 with {\arrow{>}}}]
\draw [postaction={decorate}] (cos 270,sin 270) --(2*cos 270,2*sin 270) node[midway, right]{$b$};
\end{scope}

\end{scope}

\end{scope}

\end{tikzpicture}
\caption{Depicted on the left is the attaching map of the 3-cell $B_{c^n}$ (folding the two hexagons yields a sphere). After internally collapsing $B_e$ by pushing in $d_{c^n}$, the attaching map of $B_{c^n}$ is deformed to the attaching map shown on the right.}
\label{fig:morecollapsing}
\end{figure}
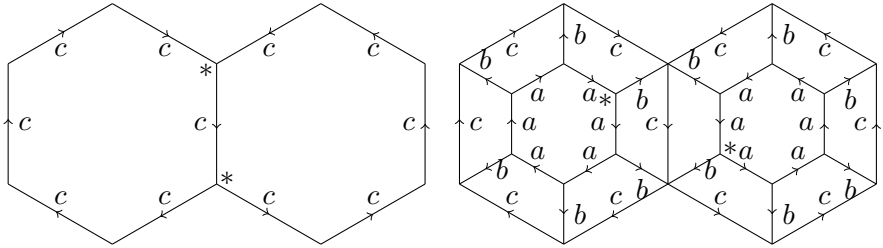

We can collapse 3-cells in $L$ following the collapsing strategy we have used on $\mathcal T$. Details on the notion of collapsing cell complexes can be found in Cohen's Book \cite{Cohen} on simple homotopy theory. Consider the 3-cell $B_e$ that is attached via the spherical diagram $S_e$. See Figure \ref{fig:3cellcollapse}. Since $c$ is a boundary vertex of $\mathcal T$, the only other 3-cell of $L$ that contains $d_{c^n}$ in its boundary is the 3-cell $B_{c^n}$ whose attaching map is the spherical power diagram associated with $c^n$. We can internally collapse $B_e$ by pushing in $d_{c^n}$. Note that $d_{c^n}$ is pushed onto $\partial B_e - d_{c^n}$, the part of the boundary that is attached via the disc diagram shown in Figure \ref{fig:S_e}. See also Figure \ref{fig:3cellcollapse}. This disc diagram consists of the power 2-cell with boundary ${a^n}$ surrounded by squares with boundary $r_e$. The attaching map of $B_{c^n}$, which is the spherical power diagram $S_{c^n}$, is deformed to the attaching map given by the spherical diagram $(S_e-\hat d_{c^n})\cup -(S_e-\hat d_{c^n})$, where $\hat d_{c^n}$ is the power 2-cell in $S_e$ with boundary word $c^n$. See Figure \ref{fig:morecollapsing}. Cancellations of squares further deforms the attaching map (but not its homotopy class) and we arrive at an attaching map that is represented by the spherical diagram $S_{a^n}$. Thus the 3-cell $B_{c^n}$ deforms into a second copy of the 3-cell $B_{a^n}$. Removing one of the $B_{a^n}$'s does not change the fundamental group or the second homotopy module of the 2-complex. 

In summary, the operation just performed removes the 3-cell $B_e$, removes the power 2-cell $d_{c^n}$, and removes the 3-cell $B_{c^n}$ without changing the fundamental group or the second homotopy module of the 3-complex. If we continue with the collapsing procedure, following the collapsing strategy for $\mathcal T$, discarding redundant power 3-cells as we go along, we end up with a 3-complex $L'$ with 2-skeleton is $K(\mathcal T, z^n)$, and with a 3-skeleton that consists of the single 3-cell $B_{z^n}$. Here $z$ is the vertex that $\mathcal T$ is collapsed to. Furthermore, $\pi_i(L')=\pi_i(L)$, $i=1,2$.  In particular, $\pi_2(L')=0$ by Theorem \ref{thm:pi2}. This gives the following result:

\begin{thm} Let $\mathcal T$ be an injective labeled oriented tree that does not contain the configurations shown in Figures \ref{fig:2-cycles}, \ref{fig:3-cycles}, \ref{fig:npc+}. Let $z$ be a vertex of $\mathcal T$. Then $\pi_2(K(\mathcal T, z^n))$ is generated by the spherical power diagram $S_{z^n}$.
\end{thm}

\section{Proof of the main Theorem}

We can now prove Theorem \ref{thm:mainthm}. Since $G({\mathcal T}, x^n)$ is not $\mathbb Z_n$ by Theorem \ref{thm:notcyclic}, it follows that $H(n)$ is not trivial. Note that $H(n)\cap \langle x \rangle=\{1\}$, since $x$ has order $n$ in the quotient $G({\mathcal T}, x^n)/H(n)=\mathbb Z_n$. Since $H(n)$ is normal it intersects all conjugates of $\langle x \rangle$ trivially. Thus, by Theorem \ref{maxfinitesubgroups}, $H(n)$ does not contain a non-trivial element of finite order. \qed

\section{Examples}

Examples of injective labeled oriented intervals ${\mathcal T}$ that do not contain edge combinations as shown in Figures \ref{fig:2-cycles}, \ref{fig:3-cycles} are contained in \cite{Rosebrock}. Hence the associated LOI-complexes $K(\mathcal T)$ are non-positively curved squared complexes. It is easy to single out examples that also avoid the edge combination show in Figure \ref{fig:npc+}. Let $\mathcal T$ be the labeled oriented tree on vertices $\{ 0, 1, ..., n-1\}$ where the edge that connects $i$ to $i+1$ is labeled with $i+3$, all numbers taken modulo $n$. Orientations of edges are irrelevant. Basic modular arithmetic shows that the labeled oriented interval $\mathcal T$ satisfies the conditions of our main theorem in case $n\ge 10$.

\begin{figure}[htbp] 
   \centering
    
\begin{tikzpicture}
\fill (0,0) circle (1pt);
\fill (1,0) circle (1pt);
\fill (2,0) circle (1pt);
\fill (3,0) circle (1pt);
\fill (4,0) circle (1pt);
\fill (5,0) circle (1pt);
\fill (6,0) circle (1pt);
\fill (7,0) circle (1pt);
\fill (8,0) circle (1pt);
\fill (9,0) circle (1pt);
\fill (10,0) circle (1pt);
\node [below] at (0,0) {$0$};
\node [below] at (1,0) {$1$};
\node [below] at (2,0) {$2$};
\node [below] at (3,0) {$3$};
\node [below] at (4,0) {$4$};
\node [below] at (5,0) {$5$};
\node [below] at (6,0) {$6$};
\node [below] at (7,0) {$7$};
\node [below] at (8,0) {$8$};
\node [below] at (9,0) {$9$};
\node [below] at (10,0) {$10$};
\draw (0,0) to (1,0);
\draw (1,0) to (2,0);
\draw (2,0) to (3,0);
\draw (3,0) to (4,0);
\draw (4,0) to (5,0);
\draw (5,0) to (6,0);
\draw (6,0) to (7,0);
\draw (7,0) to (8,0);
\draw (8,0) to (9,0);
\draw (9,0) to (10,0);
\node [above] at (0.5,0) {$3$};
\node [above] at (1.5,0) {$4$};
\node [above] at (2.5,0) {$5$};
\node [above] at (3.5,0) {$6$};
\node [above] at (4.5,0) {$7$};
\node [above] at (5.5,0) {$8$};
\node [above] at (6.5,0) {$9$};
\node [above] at (7.5,0) {$10$};
\node [above] at (8.5,0) {$0$};
\node [above] at (9.5,0) {$1$};

\end{tikzpicture}
 \caption{An example of a LOI that satisfies the condition of Theorem \ref{thm:mainthm}. Any orientation will work.}
   \label{fig:example}
\end{figure}
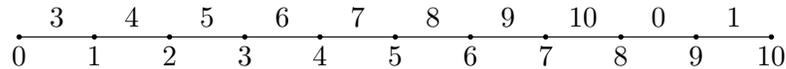

\bigskip Renata Gerecke, Pomona College, {\em renatavgerecke@gmail.com}

Jens Harlander, Boise State University, {\em jensharlander@boisestate.edu}
 
Ryan Manheimer, The College of New Jersey, {\em ryanmanheimer@gmail.com}
 
Bryan Oakley, University of Georgia,  {\em boa1540@uga.edu}

Sifat Rahman, University of Michigan-Ann Arbor,  {\em sifatj@umich.edu}

\end{document}